\documentclass[12pt]{amsart}
\usepackage{amscd,multirow,righttag,xspace}
\newcommand{\titlestring}{Crepant Terminalisations and 
   Orbifold Euler Numbers for $\SL(4)$ Singularities}
\newcommand{\runningheadstring}{Crepant Terminalisations and 
   Orbifold Euler Numbers for $\SL(4)$}
\newcommand{\msc}{(1991):  32S45 (Primary) 14L30 14E30 (Secondary)}
\newcommand{\name}{Alexander V. Sardo Infirri}

\newcommand{\ie}{i.e.\xspace}
\newcommand{\map}[5]{$$\begin{array}{rccc} 
#1\colon & #2 & \longrightarrow & #3 \\ 
         & #4 & \longmapsto & #5
\end{array}$$}

\newcommand{\Cl}{\operatorname{\rm Cl}}
\newcommand{\Cent}{\operatorname{\rm N}}

\newcommand{\Hom}{\operatorname{\rm Hom}}
\newcommand{\GL}{\operatorname{\rm GL}}

\newcommand{\SL}{\operatorname{\rm SL}}

\def\<<#1>>{\langle#1\rangle}
\def\qsing#1/#2(#3){\frac{#1}{#2}(#3)} 
\newcommand{\PP}{{\mathbb P}} 
\newcommand{\Q}{{\mathbb Q}} 
\newcommand{\R}{{\mathbb R}} 
 
\newcommand{\Z}{{\mathbb Z}} 
\newcommand{\C}{{\mathbb C}}
\newcommand{\CX}{{\C^*}}

\newcommand{\bE}{\bar E}

\newcommand{\bG}{\bar G}

\newcommand{\bV}{\bar V}

\newcommand{\bg}{\bar g}

\newcommand{\Ge}{G_\eta}
\newcommand{\He}{H_\eta}
\newcommand{\bGe}{\bar\Ge}

\newcommand{\cGo}{canonical Gorenstein orbifold\xspace}
\newcommand{\cGos}{{\cGo}s\xspace}

\newcommand{\Hs}{H${}^*$}
\newcommand{\Is}{I${}^*$}

\newcommand{\Term}{\text{Term}}

\newcommand{\wt}{\text{wt}}

\newcommand{\birationalto}{\dashrightarrow}
\newcommand{\chiorb}{\chi_{\text{orb}}}

\newtheorem{thm}{Theorem}[section]

\newtheorem{prop}[thm]{Proposition}
\newtheorem{cor}[thm]{Corollary}
\newtheorem{lemma}[thm]{Lemma}

 \theoremstyle{definition}

\newtheorem{dfn}[thm]{Definition}

\newtheorem{conj}{Conjecture}

\newtheorem{notation}[thm]{Notation}

\theoremstyle{remark}

\newtheorem{rmk}[thm]{Remark}      
\newtheorem{question}{Question} 
\newtheorem{claim}[thm]{Claim}

\begin{document}
\title[\runningheadstring]{\titlestring\footnote{{Maths Subject 
  Classification \msc}}}
\author{\name}
\email{sacha@kurims.kyoto-u.ac.jp}
\address{Research Institute for Mathematical Sciences\\ Ky\=oto University\\ 
  Oiwake-ch\protect\=o\\ Kitashirakawa\\ Saky\protect\=o-ku\\ Ky\=oto
  606-01\\ Japan}
\date{1 October 1996} 

\begin{abstract}
Let $X$ and $Y$ be two analytic canonical Gorenstein orbifolds.  A
resolution of singularities $Y\to X$ is called an \emph{Euler
resolution} if $Y$ and $X$ have the same orbifold Euler number.  If
$Y$ is only terminal rather than smooth, it is called an \emph{Euler
terminalisation}.

It is proved that Euler terminalisations exist for toric varieties in
any dimension, for 4-dimensional toroidal varieties, and for
singularities $\C^4/G$ where $G$ belongs to certain classes of
$\SL(4)$ subgroups.  The method of proof is expected to be
applicable to a sizeable number of finite $\SL(4)$ subgroups and to
lead to a generalisation of the Dixon-Harvey-Vafa-Witten orbifold
Euler number conjecture to dimension~4.
\end{abstract}
\maketitle
\tableofcontents

\setcounter{section}{-1}

\section{Introduction}
\label{sec:intro}

An analytic $n$-fold $X$ will be called a {\em canonical Gorenstein
  orbifold\/} if it has at most canonical Gorenstein singularities and
is such that, for each $x\in X$, there exists a finite group
$\pi_x<\SL(n)$ such that
$$(X,x)\cong (\C^n/\pi_{x}, 0)$$
as germs of analytic spaces.

\subsection{The Orbifold Euler Number}
\label{sec:intro:orbi}

The {\em orbifold Euler number\/} of $X$  is defined as the (finite) sum
\begin{equation}
  \label{eq:orb-eul-def}
  \chiorb(X) := \sum_{k\geq 1} k\chi(m^{-1}(k)),
\end{equation}
where $\chi$ is the ordinary Euler number and 
\map{m}{X}{\Z}{x}{|\Cl(\pi_x)|} is the upper semi-continuous map 
assigning to each point $x$ the number of conjugacy classes of 
$\pi_x$.  It is easy to show~\cite{roan:calabi-yau} that if $M$ is an 
$n$-fold admitting a $G$-action whose non-trivial elements' fixed-point 
loci have codimension at least two, and such that $M/G$ has only 
Gorenstein singularities, then
\begin{align}
\chiorb(M/G) &\phantom{:}=  \chi_{DHVW}(M;G),\\
             &:= \sum_{[g]\in\Cl(G)}\chi(M^g/\Cent^G_g)  
\end{align}
where the right-hand side denotes the Dixon-Harvey-Vafa-Witten Euler
number proposed in~\cite{dhvw:i,dhvw:ii}.

If $X$ and $Y$ are \cGos and $Y\to X$ is a bi-meromorphic map such
that $\chiorb(Y)=\chiorb(X)$, then $Y$ will be called an {\em Euler
  blow-up\/} of $X$.  If $Y$ is in addition smooth, then $Y$ will be called
an {\em Euler resolution}.  Restated in the above terminology, the
Dixon-Harvey-Vafa-Witten Euler conjecture~\cite{dhvw:i} is 
\begin{quotation}
Every $3$-dimensional \cGo has an Euler resolution.
\end{quotation}
It took ten years to give a positive answer to the conjecture
\cite{mar_ols_per,roan:mirror_cy,mark:res_168,roan:res_a5,ito:trihedral,roan:crepant}.

\subsection{Euler Terminalisations}
\label{sec:intro:euler-term}

Right from the start it was recognized that the analogous conjecture
in dimension 4 (Do all 4-dimensional \cGos have Euler resolutions?) is
trivially false: the simplest non-smooth example $\C^4/\<<-1>>$ is
already terminal.

However, rephrased slightly, the conjecture can be made to look much
more promising.  For this, note that the existence of Euler
resolutions in dimension 3 is equivalent to saying that the minimal
models for these singularities are smooth.  In other words ``smooth'' is
equivalent to ``terminal'' for 3-dimensional Gorenstein finite
quotient singularities.  

\begin{dfn}
  A Euler blow-up $Y\to X$ such that $Y$ has only terminal
  singularities is called an {\em Euler terminalisation\/} of $X$.
  The property $\Term(X)$ is defined to be true if and only if such a
  $Y$ exists.
\end{dfn}

\begin{dfn}
  The property $\Term(n)$ is defined to be true if and only if
  $\Term(X)$ is true for all $n$-dimensional \cGos $X$.
\end{dfn}

Thus, $\Term(2)$ is true in virtue of classical work and $\Term(3)$ is
true by the recent work mentioned above.   
\begin{question}
    \label{conj}
    Is $\Term(n)$ true?
\end{question}
The next open case is of course $\Term(4)$.  As in the case of
dimension~3, the question reduces to the problem of constructing Euler
terminalisations for the local singularities $\C^{4}/G$ for all the
small subgroups $G$ of $\SL(4)$.

We shall use the following terminology.  A subgroup $G < \SL(V)$ will
be called \emph{reducible} if $V$ is reducible as a $G$-module. The
{\em type\/} of a group $G<\SL(n)$ denotes the dimensions of the
irreducible representations of $G$ appearing in the chosen special
linear representation $\C^n$.  For instance, irreducible groups have
type $(n)$ and abelian groups have type $(1,1,\dots,1)$.

For any $n$, denote by $Z_n$ the cyclic central subgroup of $\SL(n)$.
For any element $g\in G$,  $\Cent^G_g$ denotes the
centralizer of $g$ in $G$, namely $\{h\in G |  h^{-1}gh=g\}$

\subsection{Main Results}
\label{sec:intro:results}

\subsubsection{Toric and Toroidal cases}

The first result is that $\Term(X)$ is true in all dimensions for
toric varieties.  Further more, it holds also for toroidal varieties
(analytic varieties which are locally isomorphic to toric varieties)
in dimension~4 and in dimension~$n$ if termination of flips can be
proved.

\begin{thm}
\label{thm:toric-toroidal}
  All simplicial toric \cGos have Euler terminalisations.
  Furthermore, if flips terminate in dimension $n$, then all
  $n$-dimensional \emph{toroidal } \cGos have Euler terminalisations.  In
  particular, this is the case in dimension 4.
\end{thm}

The proof of the toric case is straightforward; crepant blow-ups of
$X$ correspond to subdividing the first quadrant in $\R^n$ by rays
whose generators all lie in the same hyper-plane; the orbifold Euler
number of any cone is equal to its volume (meaning the volume of the
simplex spanned by its generators), and since the sum of the volumes
of the cones in the subdivision is equal to the total volume of the
original quadrant, the orbifold Euler number is seen to remain
invariant under crepant blow-ups.  The fact that among all the crepant
blow-ups there exists a terminal one is a consequence of the Toric
Minimal Model Program~\cite{reid:toricmmp}.

There is no minimal model program as yet for non-toric varieties in
dimensions 4 and over.  However, in the toroidal case, a well-known
technique makes it is possible to construct flips by patching together
local toric flips and using uniqueness.  Thus if termination of flips
is also known, (as it is in dimension 4~\cite{kmm:intro_mmp}), the
existence of the terminal model follows and this again must have the
same orbifold Euler number as the original variety.

\subsubsection{$\SL(4)$ subgroups of type $(3,1)$}

The second set of results concerns 4-dimensional non-abelian
singularities created by finite $\SL(4)$ subgroups of type (3,1).  

Let $G<\SL(4)$ be a finite subgroup which stabilises a line
$V^2\subset V=\C^{4}$ and let $V^1$ be a $G$-submodule such that
$V=V^1\oplus V^2$.  Denote by $\eta$ the generic point of the line
$\{0\}\times V^2$, and by $\Ge$ the stabiliser of $\eta$.  Note that
$\Ge$ is a subgroup of $\SL(V^1)\times\{1\}\cong \SL(3)$.

\begin{thm}
  If $G<\SL(4)$ fixes a line and is such that the group $\Ge$ does not
  contain $Z_3$ as a subgroup, then $\C^4/G$ has an Euler
  terminalisation with only toric singularities.
\end{thm}

The method of proof essentially consists in using the results of
Roan \cite{roan:calabi-yau} to construct a 3-dimensional Euler
resolution of $\C^3/\Ge$ which is equivariant under the larger group
$G$. 

As mentioned in the next section, the assumption that $\Ge$ does 
not contain the group $Z_3$ is not essential to the method.  However, 
so far the author has been unable to construct the equivariant 
resolutions without it.  An attempt is made in Section~\ref{sec:3-1:centre:2}.

\subsubsection{$\SL(4)$ subgroups containing $Z_4$}

The third set of results reduces the question $\Term(\C^4/G)$ for the
irreducible $G$ which contain $Z_4$ to a conjecture regarding the
existence of the equivariant resolutions of $\SL(3)$ mentioned above.
The conjecture is proved for irreducible subgroups of $\SL(3)$ which
do not contain $Z_3$, but remains open in the general case.

\begin{rmk}
The material presented here can no doubt be pushed further, but the
author has so far been unable to do so.  Nevertheless, it is hoped
that the attentive reader will be able to perceive a direction in
which to proceed.  Is is also conceivable that the general strategy
that emerges from this approach may be applicable to an understanding
of quotient singularities in higher dimensions.
\end{rmk}

The idea here is as follows.  Suppose that $G<\SL(n)$ acts on $\C^{n}$
and contains the centre $Z_{n}$, and write $\bG:=G/Z_n$ and $\bV :=
\text{Bl}_0V/Z_n$.  Then Lemma~\ref{lemma:bV} implies that $\bV/\bG$
is another \cGo and $\bV/\bG \to V/G$ is an Euler blow-up.  With the
aid of another Lemma (The Patching Lemma~\ref{lemma:patching}), the
problem is thus reduced to the construction and patching of Euler
blow-ups of local neighbourhoods of $\bV/\bG$.  The advantage of these
is that the stabilisers of $\bG$ in the tangent space to the blowup
$\bV$ are simpler than those of $G$ (because the stabiliser of $\bG$
must fix a line in the tangent space to $\bV$, so its type must be
$(t,1)$, where $t$ is a partition of $n-1$).
  
This approach is spelt out in Section~\ref{sec:4-centre} for
irreducible subgroups of $\SL(4)$ which contain $Z_{4}$; it reduces
the problem to constructing and patching together equivariant $\SL(3)$
resolutions.

\subsubsection{Discussion}

What if, on the other hand, the group $G$ is irreducible, but doesn't
contain the centre of $\SL(V)$?  

A complete answer to Question~\ref{conj} for all finite $\SL(4)$
singularities seems out of reach of the methods suggested here, if only
because the cases when $G$ does not contain $Z_4$ include simple
groups, such as the alternating group $A_5$; if dimension~3 is any
indication~\cite{mark:res_168,roan:res_a5}, it seems that ad-hoc
methods will be necessary to construct a terminalisation.

However, the possibility is open that, as a general rule, the
non-simple finite subgroups $\SL(n)$ not containing $Z_n$ are few in
number and relatively amenable in form.  For instance, if $\dim V=2$
no cases occur.  For $\dim V=3$, one has --- apart from the simple
groups (H) and (I) which require ad-hoc methods --- ``half'' the
groups of type (C) and ``half'' of those of type (D).  It turns out
that under the assumption $Z_3\not < G$ these are semi-direct products
of abelian groups with the alternating group $A_3$ and the symmetric
group $S_3$ respectively.  This allows one to construct their Euler
resolutions from toric resolutions --- see~\cite{roan:calabi-yau}.
Extension of this method to the $\SL(4)$ case would seem to be
feasible.
  
Table~\ref{tab:sl4} outlines the state of the $\Term(\C^{4}/G)$ 
problem. 

\begin{table}[htbp]
  \begin{center}
    \leavevmode
\setlength{\fboxrule}{0pt}
\begin{tabular}{|l|l|l|l|l|}
        \hline
        \multicolumn{3}{|c|}{Type of Group} & Method & Status\\
        \hline
        \multirow{7}{12mm}{Irred.\ \textbf{(4)}} &
        \multirow{2}{16mm}{$G\not > Z_{4}$} 
        & \multirow{1}{22mm}{Simple} & \fbox{Ad hoc method?} &\\
        \cline{3-5}
        & & \multirow{1}{22mm}{Non-Simple} & 
        \fbox{\parbox{3.2cm}{Small \# of cases?  
            Semi-direct products of abelian groups with $A_{4}$ and 
            $S_{4}$?}} &  \\
        \cline{2-5}
        &  \multicolumn{2}{l|}{$G > Z_{4}$}  & 
        \fbox{\parbox{3.2cm}{4-d Euler blow-up of origin. Patch
        together equivariant lower-d blowups.}} &
        \S\ref{sec:4-centre}\ \ \ ?\\
        \hline
        \multirow{5}{12mm}{Red.} & \multirow{2}{16mm}{\textbf{(3,1)}}
        & \multirow{1}{22mm}{$ \Ge\not >Z_{3}$} &
        \fbox{\parbox{3.2cm}{Construct equivariant 3-d Euler
        resolution. Reduces to the toric case.}} 
        & \S\ref{sec:3-1:nocentre} OK\\
        \cline{3-5}
        & & \multirow{1}{22mm}{$\Ge >Z_{3}$} & \fbox{\parbox{3.2cm}{Euler 
        blow-up of fixed line. Reduce to 2-d equivariant blow-up.}} &  
        \S\ref{sec:3-1}\ \ \ ?\\
        \cline{2-5}
        & \multicolumn{2}{l|}{\textbf{(2,2)}}  & \fbox{Use 2-d results} &  
        \S\ref{sec:2-2} ?\\
        \cline{2-5}
        & \multicolumn{2}{l|}{\textbf{(2,1,1)}}   
        & \fbox{Use 2\&3-d results} & \S\ref{sec:2-1-1} ?\\
        \cline{2-5}
        & \multicolumn{2}{l|}{\textbf{(1,1,1,1)} (Abelian)}  & \fbox{Toric 
        MMP} & \S\ref{sec:toric}\ \ \ OK\\
        \hline
\end{tabular}
\vspace{1em}
\caption{Constructing Euler Terminalisations of $\C^4/G$ 
for $G<\SL(4)$.  (Bold numbers in parentheses indicate the {\em 
type\/} of the group.  In the cases $(3,1)$, the group $\Ge$ denotes 
the stabiliser in $G$ of a generic point of the line fixed by $G$.  )}
    \label{tab:sl4}
  \end{center}
\end{table}
The {\em Status\/} column indicates the section of this paper which
makes some contribution to the problem.  Questions which are solved in
this paper are indicated by the mention ``OK''.  A question mark
indicates that the question is still open and that methods of this
paper are relevant.  Finally, where no results are known, nothing
is indicated in the Status column, but some guesses are given as to
the likely situation, based on the present state of knowledge.

\subsection{Open Problems}
\label{sec:open}

Interesting open problems arise in relation to recent work of Ito and
Reid~\cite{reid_ito} which establishes a one-one correspondence
between crepant divisors of $V/G$ and conjugacy classes ``of weight
1'' in (the dual group to) $G$.

One interesting question is how this correspondence behaves under the
Euler blow-ups which are constructed for the $\SL(4)$ singularities
mentioned above.  A deeper understanding of this (apart from being of
interest in itself) would also no doubt allow one to say more about
the singularities of the terminalisation.

For instance, for toric and toroidal varieties, the Euler blow-up can
be made projective rather than just analytic.  But in general, the
gluing process in the Patching Lemma~\ref{lemma:patching} does not in
itself guarantee that the blow-up $X$ will be projective, because a
divisor which is reducible when considered locally in the
neighbourhood of one of the points $[\xi]$ might have two of its
components identified by the gluing process.  Such troublesome cases
might conceivably be ruled out by a deeper understanding of the
correspondence between divisors and conjugacy classes.

\subsection{Outline}
\label{sec:intro:outline}

Section~\ref{sec:toric} deals with the construction of Euler
terminalisations for the toric and toroidal cases.
Section~\ref{sec:3-1} deals with the finite $\SL(4)$ groups of type
(3,1).  Section~\ref{sec:blowing-up} presents a blowing-up
construction and some other technical lemmas which are then used in
Section~\ref{sec:4-centre} in dealing with irreducible finite $\SL(4)$
subgroups containing $Z_4$.  Finally, Section~\ref{sec:2-} makes some
comments regarding the finite $\SL(4)$ groups of type (2,2) and
(2,1,1).

\subsection{Acknowledgments}
\label{sec:intro:ack}

I wish to acknowledge S.Mori and S-S.Roan for the many ideas and
suggestions they contributed to this paper.

This research was undertaken at the Research Institute for
Mathematical Sciences of Kyoto University thanks to an European
Commission Science and Technology Fellowship; I am grateful to my host
institution and its staff for their hospitality and to the European
Commission for their financial support.

\section{Toric and Toroidal Cases}
\label{sec:toric}

\begin{thm}[Toric Minimal Model Program]
\label{thm:toric-mmp}
Let $Y$ be a toric variety and $X $ be a simplicial toric variety
which admits a projective birational toric morphism $f\colon X \to Y$.
Then there exists a sequence $X \stackrel h\birationalto
Z\xrightarrow{g}Y$ such that
  \begin{enumerate}
  \item $h$ is a composite of toric divisorial contractions or toric flips.
  \item $g$ is a projective morphism and $Z$ is a simplicial toric
    variety with terminal singularities such that $K_Z$ is relatively
    nef for $g$.
  \end{enumerate}
  Note that if $Y$ has canonical singularities, then $K_Z=g^*K_Y$ in
  the sense of $\Q$-Cartier divisors, \ie  $g$ is crepant.
\end{thm}
\begin{proof}
  See~\cite[Theorem 0.2]{reid:toricmmp}, where the result is proved
  under the assumption that $Y$ is projective.  As
  remarked in~\cite{reid:toricmmp} this assumption is not essential and
  the result is valid for non-complete toric varieties also (the
  easiest way to see this is to reduce the non-projective case to the
  projective one by completing the fan in an appropriate way).
\end{proof}

\begin{cor}
\label{cor:abelian-terminal}
 All toric \cGos admit (toric) Euler terminalisations.
\end{cor}
\begin{proof}
  Note that a toric variety has at most orbifold singularities if and
  only if it is simplicial.  It is therefore sufficient to prove that
  any crepant blow-up of a simplicial Gorenstein toric variety must
  have the same orbifold Euler number as the original.  But this is
  true because the orbifold Euler number of a simplicial toric variety
  is just the  volume\footnote{Also called the {\em multiplicity\/} of 
  the cone.} of the cone, meaning the volume
  of the simplex defined by the cone's generators; a crepant blow-up
  corresponds to a fan subdivision by one-dimensional rays whose
  primitive generators all belong to the same plane, and therefore the
  sum of the volumes of the cones in such a subdivision is equal
  to the volume of the original cone.
\end{proof}

The minimal model program for general varieties is at present only
proved in dimension~3.  In dimension~4, although termination has been
shown~\cite{kmm:intro_mmp} existence of flips remains a problem.
Nevertheless, a technique well-known to minimal model program
specialists allows one to use the Theorem above and the termination
result to say something about toroidal varieties, \ie  varieties
which are only locally isomorphic to toric varieties.  The argument
can phrased for general $n$, even though at present, termination has
only been proved for $n\leq 4$.

\begin{thm}
\label{thm:toroidal-mmp}
Assume that flips terminate in dimension $n$.  Then all
$n$-dimensional toroidal \cGo admit Euler terminalisations (which are
themselves toroidal).  In particular, this is true in dimension~4.
\end{thm}
\begin{proof}
  Let $Y$ be a \cGo locally isomorphic to a toric variety and let
  $p\colon X\to Y$ be any resolution obtained by toric blowups.
  
  Suppose that $K_X$ is not $p$-nef, and let $c\colon X \to W$ be an
  extremal contraction.  If it is a divisorial contraction, then
  replace $X$ by $W$.  If $c$ is a small contraction, consider its
  restriction $C_{U_X}\colon U_X \to U_W$ to the inverse image of a
  local toric neighbourhood $U_Y\subset Y$.
  Theorem~\ref{thm:toric-mmp} implies the existence of a local flip
  $c^+_{U_X}\colon U^+_X \to U_W$ over each neighbourhood $U_W$.  A
  flip being unique if it exists, the local flips patch together on
  the overlaps to form a global flip $c^+\colon X^+\to W$.  Thus,
  existence of flips is guaranteed in this case.  Applying the same
  procedure repeatedly (and using the termination hypothesis) results
  in a projective morphism $p\colon Z\to Y$ such that $Z$ has
  $\Q$-factorial terminal singularities with $K_Z$ being $p$-nef,
  which means that $p$ is crepant, since $Y$ is canonical.
  Furthermore, any toric terminalisation which is crepant must have
  the same orbifold Euler number by the volume argument in
  Corollary~\ref{cor:abelian-terminal}.
\end{proof}

Corollary~\ref{cor:abelian-terminal} and
Theorem~\ref{thm:toroidal-mmp} together give 
Theorem~\ref{thm:toric-toroidal}.

\section{Groups $G$ of type (3,1)}
\label{sec:3-1}

\subsection{Notation}
\label{sec:3-1:notation}

Let $G<\SL(n+1)$ be a finite subgroup such that $V=\C^{n+1}$
decomposes into two irreducible $G$-modules: $V=V^1\oplus V^2$, with
$V^1$ of dimension $n$, and $V^2$ of dimension 1.

Denote by $\eta$ the generic point of the line $\{0\}\times V^2$, and
by $\Ge$ the stabiliser of $\eta$.  Note that $\Ge$ is a subgroup of
$\SL(n)\times\{1\}\cong \SL(n)$.  The quotient $C:= G/\Ge$ is cyclic,
being a naturally a subgroup of $\GL(V^2)\cong\CX$.  The canonical
quotient map is denoted $\pi\colon G\to C$ and the induced map on
conjugacy classes is denoted $\pi_*\colon\Cl(G)\to
\Cl(C)$.

Note that $G$ can be considered as a subgroup of $\GL(V^1)=\GL(n)$ by
forgetting about the last row and column of any matrix element.  Let
$h\in G$ be an element such that $G=\<<\Ge,h>>$ --- \ie a
representative for a generator of $C$ --- and denote by $h_1$ the
$n\times n$ sub-matrix consisting of the first $n$ rows and columns of
$h$.  Let $\lambda_h$ be a complex $n$-th root of $\det(h_1)^{-1}$.
Then $h':= \lambda_h h_1\in\SL(n)$ and normalises $\Ge$.  Hence $\Ge$
is a normal subgroup of $G':=\<<\Ge,h'>> < \SL(n)$ with quotient a
cyclic subgroup $C'$.  Note that $G'$ and $C'$ are defined up to a
choice of the $n$-th root $\lambda_h$.

The canonical quotient map is denoted $\pi'\colon G'\to C'$ and the
induced map on conjugacy classes is denoted
$\pi^{\prime}_*\colon\Cl(G')\to \Cl(C')$.

\subsection{Equivariant Resolution}
\label{sec:3-1:resolution}

From now on, let $n=3$, \ie  consider the case where $G<\SL(4)$ fixes
a line in $\C^4$ (and therefore $\Ge<\SL(3)$).

\begin{conj}
  \label{conj:3-1-equiv-Euler}
  Let $G<\SL(4)$ be a finite subgroup which stabilises a line
  $V^2\subset V=\C^4$ and let $\Ge$ denote the stabiliser in $G$ of
  the generic point of $V^{2}$ and let $G'$ and $C'$ be defined as in
  Section~\ref{sec:3-1:notation} above.
  
  Then there exists a $G'$-invariant Euler resolution $W^1\to
  V^1/\Ge$, satisfying
  \begin{equation}
    \label{eq:chi-pi}
    \chi((W^1)^{c'}/{C'})=|\pi^{\prime -1}_*([c'])|,
  \end{equation}
  for all $c'\in C'$.
\end{conj}

If true, we shall see that this conjecture implies that
$\Term(\C^{4}/G)$ is true for groups stabilizing a line (\ie  types
$(3,1)$ and $(2,1,1)$) and (cf.\ Section~\ref{sec:4-centre}) for all
irreducible $G$ which contain $Z_{4}$.

\begin{prop}
\label{prop:3-1-invariant-resolution-gives-terminalisation}
Suppose that Conjecture~\ref {conj:3-1-equiv-Euler} is true for $G$.   
Then $W^1\times V^2/C\to V/G$ is an Euler blowup, 
and $V/G$ has an Euler terminalisation with only
toric singularities.
\end{prop}
\begin{proof}
    First, note that since $V^1/\Ge$ has dimension at most three, it has
  a minimal model.  The fact that there are only a finite number of
  distinct minimal models implies that $\CX$ acts on any of them: for
  if the action of some element $\lambda\in\CX$ produced a different
  minimal model, then, by continuity, one could produce countably many
  distinct minimal models by acting with a countable family of
  distinct neighbours of $\lambda$. 
  
  Second, note that the equality~\eqref{eq:chi-pi} would follow from
  the same equality with $C'$ replaced by the group $C$.  For if
  $\phi$ is any element of $G$, one can find $\lambda\in\CX$ such that
  $\phi'=\lambda \phi\in G'$.  The invariant sets $(W^1)^\phi$ and
  $(W^1)^{\phi'}$ have the same homotopy type, by an easy application
  of Bialynicki-Birula's well-known decomposition
  theorem~\cite[Thm.~4.1]{bb:algebraic_groups} to the smooth variety
  $W^1$.  Hence $\chi((W^1)^\phi)=\chi((W^1)^{\phi'})$ and so,
  averaging,
  \begin{equation}
  \label{eq:chi-equality}
  \chi((W^1)^\phi/C)=\chi((W^1)^{\phi'}/C').
  \end{equation}
  On the other hand, if $c$ and $c'$ denote the images in $C$ and $C'$ 
  of $\phi$ and $\phi'$ respectively,
  \begin{equation}
  \label{eq:pi-equality}
  |\pi_*^{-1}([c])|=|\pi^{\prime -1}_*([c'])|,
  \end{equation}
  where $\pi$ denotes the projection $\pi\colon G\to C$ and
  $\pi_*\colon\Cl(G)\to\Cl(C)$ the induced map on
  conjugacy classes, and similarly for the primed symbols.
  
  Hence one can assume formula~\eqref{eq:chi-pi} to be valid for the
  group $C$.  The variety $(W^1\times V^2)/G$ is a blow-up of $V/G$
  having only cyclic quotient singularities resulting from the
  residual action of $C=G/\Ge$ on $W^1\times V^2$.

 Its  orbifold Euler number can be expressed as a sum:
  \begin{align}
    \chiorb(W^1\times V^2/C) &=\sum_{[c]\in\Cl(C)}
    \chi((W^1\times V^2)^c/C),\quad\text{since $\Cent^C_c=C$}\notag\\ 
               &= \sum_{[c]\in\Cl(C)} \chi((W^1)^c/C),\\
               &\qquad\text{since $(V^2)^c$ is contractible}. 
    \label{eq:orb-sum}
  \end{align}
  
  On the other hand, one has:
   \begin{equation}
    \label{eq:cclass-sum}
   |\Cl(G)|=  \sum_{[c]\in \Cl(C)}  |\pi_*^{-1}([c])|, 
  \end{equation}
  which agrees with the previous sum term-by-term.  Hence, $Y:=(W^1\times
  V^2)/G\to (V^1\times V^2)/G$ is an Euler blow-up with only toric
  (cyclic) singularities.
  
  Applying the minimal model program (Theorem~\ref{thm:toroidal-mmp})
  to $Y$, one obtains a crepant terminalisation $t\colon Z\to Y$ which
  satisfies $\chiorb(Z)=\chiorb(Y)=\chiorb(V/G)$, and has only toric
  singularities.
\end{proof}

Thus in the case where $G$ fixes a line in $V=\C^4$ it suffices to
prove the existence of a $G'$-equivariant Euler resolution $W^1\to
V^1/\Ge$ which satisfies equation~\eqref{eq:chi-pi}.

\begin{rmk}
  The same method as above can be used to deal with the easier case
  when $n=2$ and $G<\SL(3)$ and fixes a line in $\C^3$.  In 2
  dimensions, the minimal model is unique, so there is no need to
  check $G'$-stability of the Euler resolution of $V^1/G_\eta$.
\end{rmk}

\subsection{Case where $\Ge$ doesn't contain $Z_3$}
\label{sec:3-1:nocentre}

In this section, Conjecture~\ref{conj:3-1-equiv-Euler} is proved in
the case where $\Ge$ is irreducible and does not contain $Z_{3}$.

In order to construct a $G$-equivariant Euler resolution $W^1\to
V^1/\Ge$, the cases to be considered are first restricted using the
following lemma (which makes use of the classification of small finite
sub-groups of $\SL(3)$ --- see~\cite{yau_yu}, although the notation
adopted here is that of~\cite{roan:calabi-yau}, which is slightly
different).

\begin{lemma}
\label{lemma:3-1-nocentre-exact}
 Suppose that 
  \begin{equation}
    \label{eq:3-1-nocentre-exact}
    1 \to \Ge \to G' \xrightarrow{\pi'} C' \to 1,
  \end{equation}
  is an exact sequence of finite groups of $\SL(3)$ 
  such that $C'$ is cyclic and non-trivial, and $\Ge$ is non-abelian and
  doesn't contain $Z_3=\<<\omega_3>>$. 
 Then exactly one of the following is true.
 \begin{enumerate}
 \item $\Ge$ is of  type (B) and $G'$ is of type (B).
 \item $\Ge$ is of  type (C), (D), (H) or (I) and  $G' =
   \<<\Ge,\omega_3>>$.
 \item $\Ge$ is of type (C) and $G'$ is of type (D), with
   $G'=\<<\Ge,R>>$ or $G'=\<<\Ge,\omega_3R>>$.
 \end{enumerate}
\end{lemma}
\begin{proof}
  The groups of type (E), (F), (G), (\Hs), (\Is) all contain $Z_3$, so
  do not occur as the group $\Ge$ by assumption.  The only finite subgroup of
  $\SL(3)$ containing the simple group (H) (resp.\ the simple group
  (I)) as a normal subgroup is (\Hs) (resp.\ 
  (\Is))~\cite[p.36]{yau_yu}.  Thus for these, the result follows
  immediately. 

  If $\Ge$ has type (B), a simple argument~\cite[\S1.4, p.18]{yau_yu}
  shows that $G'$ must also have type (B).
  
  It remains to deal with the case where $\Ge$ has type (C) or (D).
  Throughout the rest of this proof, write
$$T:=
\begin{pmatrix}
0& 1& 0\cr 0& 0& 1\cr 1& 0& 0
\end{pmatrix},
$$ for the element which, together with a diagonal group, generates a
group of type (C).  To get a group of type (D), recall that one must
add to a group of type (C) an element of the form 
  \begin{equation}
    \label{eq:phi}
    \phi= \begin{pmatrix} a& 0& 0\cr 0& 0& c\cr 0& b& 0\cr
            \end{pmatrix} \quad\text{with }abc=-1. 
  \end{equation}

  \begin{claim}
    If $\Ge$ has type (C) or (D), then $G'$ must also be of type (C)
    or (D) (though not necessarily the same type as $\Ge$).
  \end{claim}
  \begin{proof}
    Denote by $x_1, x_2, x_3$ the standard coordinates on $\C^3$. If
    $\Ge$ is of type (C) or (D) and does not contain the centre $Z_3$,
    then the monomial $x_1x_2x_3$ is invariant under $\Ge$ up to
    scale~\cite[\S 1.3]{yau_yu}.  It follows from an easy argument
    that $G'$ must also leave $\C x_1x_2x_3$ invariant.  Thus $G'$
    cannot be primitive.  This means that $G'$ must be of type (C) or
    (D).
  \end{proof}

  The next step is to study the normal diagonal subgroups of $\Ge$ and
  $G'$, which are denoted by $\He$ and $H'$ respectively. 

  \begin{notation}
    The standard toric notation for diagonal matrices will be used:
$$\qsing1/d(r_1, r_2, \dots, r_n):\equiv [\exp({\frac{2\pi i r_1}{ d}}),
\exp({\frac{2\pi i r_2}{d}}), \dots, \exp({\frac{2\pi i r_n}{d}})]$$
  \end{notation}

  \begin{claim}
    \label{claim:order3}
    If $H'$ contains an element of order $3$ then that element must be
    $\omega_3$ or $\omega_3^2$.  As a consequence, all the elements of
    $\He$ have orders prime to $3$.
  \end{claim}
  \begin{proof}
    If $x\in H'$ has order $3$ and its does not belong to the centre
    $Z_3$ then it can be chosen to be of the form
    $x=\qsing1/3(i,i+1,i+2)$ for some $i\in\{0,1,2\}$. But then
  \begin{align*}
 (TxT^{-1})x^{-1} &= \qsing1/3(i+1,i+2,i) - \qsing1/3(i,i+1,i+2)\\
                  &= \qsing1/3(1,1,-2)\\
                  &= \qsing1/3(1,1,1)=\omega_3
  \end{align*}
  On the other hand, since $x$ normalises $\Ge$, one has 
  $xT^{-1}x^{-1}\in\Ge$ and so $\omega_3=TxT^{-1}x^{-1}\in\Ge$, which
  contradicts the hypothesis of the lemma.  Thus the only elements of
  order $3$ in $H'$ are $\omega_3$ and $\omega_3^2$. As a consequence,
  $\He$ has no elements of order $3$, since $Z_3\not < \He$.  The
  claim for $\He$ follows immediately from this.
  \end{proof}

  \begin{claim}
  \label{claim:hprime}
    $H' < \<< \He,\omega_3 >>$, i.e. $H'$ is either equal to
    $\He$ or equal to $\<< \He,\omega_3>>$.
  \end{claim}
  \begin{proof}
    Let $\varphi\in H'$.  I begin by showing that $\varphi^3\in \He$.

    Since $H'$ is normal in $G'$, the element $T \varphi T^{-1}$ is
    diagonal, and so, therefore, is $f:=\varphi^{-1}T \varphi T^{-1}$.
    One has $fT = \varphi^{-1} T\varphi\in \Ge$, since $\Ge$ is
    normal in $G'$, so $f$ belongs to $\Ge$.  Since $f$ is also
    diagonal, it follows that $f\in\He$.

      Writing $\varphi = \qsing 1/d(a,b,-a-b)$, one has
      \begin{align*}
        f &=\varphi^{-1}(T \varphi T^{-1})=
        \qsing1/d(-a,-b,a+b)+\qsing1/d(b,-a-b,a)\\ 
        &=\qsing1/d(b-a,-a-2b,2a+b)\\ \intertext{and} 
    T^{-1}f T &=      \qsing1/d(2a+b,b-a,-a-2b).
      \end{align*}
      Dividing the second element by the first gives
  $$T^{-1} f T f^{-1} = \qsing1/d(3a,3b,-3(a+b))= \varphi^3,$$
  so $\varphi^3\in\He$.

  Let $x:=\varphi^{-3}\in\He$. By Claim~\ref{claim:order3}, the
  order of $x$ is prime to 3, so there exists an integer $l$ such that
  $x^{3l}=x$.  Writing $\alpha := x^l=(\varphi^{-3l})\in\He$, one has
  $(\alpha\varphi)^3=1$, and so Claim~\ref{claim:order3} again implies
  that $\alpha \varphi =\omega_3$ or $\omega_3^2$.  
  Thus $H'=\<<\He,\omega_3>>$.
  \end{proof}

  Now one can deal with the groups $\Ge$ and $G'$ themselves.
  To begin with, since the quotient $C'=G'/\Ge$ is cyclic, $G'=\<<
  \Ge,\phi>>$ for some $\phi\in G'\setminus \Ge$.  Now any
  element $\phi$ in a group of type (C) or (D) has associated to it a
  permutation $\sigma(\phi)\in S_3$, defined according to how it
  permutes the coordinates $x_1$, $x_2$, $x_3$.  If $\sigma(\phi)$ is
  the identity, then $\phi$ is diagonal, whereas if $\sigma(\phi)$ is
  a permutation of order 3 then $\phi T$ or $\phi T^{-1}$ is diagonal.
  In these cases, since $T\in\Ge$, it follows from the claim above
  that $G' = \<< \Ge,\omega_3>>$.

  The only remaining possibility is that $\sigma(\phi)$ equals a
  transposition or order 2, which  can be assumed to be the
  transposition $(12)$, by
  multiplying $\phi$ by a suitable power of $T$.  Thus $\phi$ is of
  the form~\eqref{eq:phi}.

  \begin{claim}
    \label{claim:type_D}
    For any $\phi$ of the form~\eqref{eq:phi}, define 
$$\tilde\phi:=
    T^{-1}\phi^2T\phi=
    \begin{pmatrix}
      A& 0& 0\cr 0& 0& C\cr 0& B& 0\cr
    \end{pmatrix},$$ with $A=-1$, $B=b^2c$ and $C=-B^{-1}$. 
    Suppose that  $Z_3\not<\<<\tilde\phi,T>>$.  Then 
    there exists an element $t\in \<<T,\tilde\phi>>$ such that 
    $t\tilde\phi=R$, where
    $$R:=\begin{pmatrix}-1& 0& 0\cr 0& 0& -1\cr 0& -1& 0
    \end{pmatrix}.$$
  \end{claim}
  \begin{proof}
   Define $f:=\tilde\phi
    T\tilde\phi^{-1}T=[-B,-B,B^{-2}]$ and
    $f':=fT^{-1}fT=[1,-B^3,-B^{-3}]$.  If the order of $B$ is a
    multiple of $3$, say $m=3k$, then $f^{2k}=\omega_3$ or~$\omega_3^2$,
    so $Z_3<\<<\tilde\phi,T>>$.  Thus if
    $Z_3\not<\<<\tilde\phi,T>>$ then the order of $B$ is prime to
    three, and a suitable power of $f'$ gives the required element
    $t=[-1,-B^{-1},-B]$, which satisfies  $t\tilde\phi=R$.
  \end{proof}
  Now the element $\phi^2$ is diagonal, so belongs to $\<< \He,
  \omega_3>>$, by claim~\ref{claim:hprime}.  Note also that
  $\phi=\tilde\phi(T\phi^{-2}T^{-1})\in\tilde\phi H'$.

{\bf Case 1: $\phi^2 \in \He$}
    In this case, $C'$ has order 2, so $\omega_3\not\in H'$ --- \ie  
    $H'=\He$.  By Claim~\ref{claim:type_D}, $t\tilde\phi=R$ for some
    element $t\in \<<T,\tilde\phi>> < \Ge$.  Thus
    $G'=\<<\Ge,\phi>>=\<<\Ge,\tilde\phi,H'>> = \<<\Ge,R>>$, since $H'=\He$.

{\bf Case 2: $\phi^2\not\in\He$}
    Then $\phi^2=\varphi\omega_3^k$ for $k=1$ or~$2$ and
    $\varphi\in\He$.  Let $\phi':=\omega_3^k\phi$.  Then
    $(\phi')^2=\varphi^2\in\He$, so by the preceding case, one may
    assume $t\tilde\phi'=R$ for some
    $t\in\<<T,\tilde\phi'>>< G_\eta$. Now
    $\tilde\phi=\widetilde{\phi'}$, 
    so $G'=\<<\Ge,\phi>>=\<<\Ge,\tilde\phi,H'>> = \<<\Ge, Rt^{-1},H'>> 
    = \<<\Ge,R,\omega_3>>$.

  This completes the proof of the lemma.
\end{proof}

\begin{prop}
\label{prop:3-1-nocentre-equiv-Euler}
Conjecture~\ref{conj:3-1-equiv-Euler} is true when $\Ge$ is
irreducible and doesn't contain $Z_{3}$.
\end{prop}
\begin{proof}
  Recall from the proof of
  Prop.~\ref{prop:3-1-invariant-resolution-gives-terminalisation},
  that $\CX$ acts on any Euler resolution.  Thus, in all cases where
  $G'=\<<\Ge,\omega_3>>$, any smooth crepant resolution of $V^1/\Ge$
  admits a $G'$-action --- indeed a $G$-action, as remarked in the
  proof of Prop.~\ref{prop:3-1-invariant-resolution-gives-terminalisation},
  since $G<\<<G',\CX>>$.
  
  Hence, Lemma~\ref{lemma:3-1-nocentre-exact} implies that it suffices
  to deal with the cases where $G=\<<\Ge,\omega_3>>$ or
  $\Ge=\<<\He,T>>$ is of type (C) with $G'$ of type (D), either equal
  to $\<<\Ge,R>>$ or equal to $\<<\Ge, R,\omega_3>>$.

  In these cases (see~\cite{roan:calabi-yau}) an Euler resolution of
  $\C^3/\Ge$ is obtained by taking a toric Euler resolution of
  $\widetilde{\C^3/\He} \to \C^3/\He$ which is $\<<T>>$-stable and
  then resolving the singularities of $\widetilde{\C^3/\He}/\<<T>>$.
  The existence of a $T$-stable toric Euler resolution follows from
  the fact that $\He$ is normal in $\Ge$.  However, $\He$ is also
  normal in $G'$, so $\widetilde{\C^3/\He}$ can also be chosen to be
  $R$-stable.  The singularities of $\widetilde{\C^3/\He}/\<<T>>$ are
  fixed points of $R$, so resolving them gives the desired
  $G'$-invariant resolution of $\C^3/\Ge$.

  The proof that $\chi((W^1)^{c'}/C')=|\pi_*^{\prime -1}([c'])|$ is done 
  by  treating case by case the three possibilities for $G'$ given by
  Lemma~\ref{lemma:3-1-nocentre-exact}.
  
  {\bf Case 1: $G'=\<<\Ge,\omega_3>>$.} In this case, since
  $\omega_3\in\CX$, the same argument as the second paragraph of the
  proof of
  Proposition~\ref{prop:3-1-invariant-resolution-gives-terminalisation}
  gives $\chi((W^1)^{\omega_3})=\chi(W^1)=|\Cl(\Ge)|$.  On the other
  hand, $G'$ is just a direct product of $\Ge$ and $Z_3$, so $\pi_*'$
  is everywhere $|\Cl(\Ge)|:1$.
  
  {\bf Case 2: $G'=\<<\Ge,R>>$.} If we denote by $\He$ and $H'$ the
  normal diagonal subgroups of $\Ge$ and $G'$ respectively, then they
  are equal.  Since they do not contain $\omega_3$, their order is
  $d^2$ (for some $d$ prime to~3) and they are a semi-direct factor in
  $G'$~\cite[Lemma 10]{roan:crepant}.  The inverse image of the
  trivial class in $C'$ is the number of $G'$-conjugate elements in
  $\Ge$.  Since $\Ge$ is normal, this is the same as the number of
  $\Ge$-conjugate elements in $\Ge$, \ie  equal to
  $\Cl(\Ge)=\chi(W^1)$.  For the non-trivial class $[R]$,
  \cite[Formula~(32)]{roan:crepant} implies that $\chi((W^1)^R/\<<R>>)=d$ and
  the proof of ~\cite[Lemma 10]{roan:crepant} again gives
  $\pi_*^{\prime-1}([R])=|Z_R\cap H'|=d$.
  
  {\bf Case 3: $G'=\<<\Ge, \omega_3 R>>$.} As we remarked in Case~1,
  the Euler number does not depend on scalar factors, so
  $\chi((W^1)^{\omega_3 R})=\chi((W^1)^R)$.  On the other hand, the
  discussion in the second paragraph of this proof implies that
  $|\pi^{\prime-1}_*([\omega_3 R])|=|\pi^{\prime-1}_*([R])|$. The
  result thus follows from Case~2.
\end{proof}

Hence, if $\Ge$ is irreducible and does not contain $Z_3$, $V/G$ has
an Euler terminalisation with only toric singularities.

\begin{question}
  The proof of Propositions~\ref{prop:3-1-nocentre-equiv-Euler} 
  depends on the classification of $\SL(3)$ groups.  Is it possible to 
  find a proof which doesn't depend on the classification?
\end{question}

\subsection{Case where $\Ge$ contains $Z_3$}
\label{sec:3-1:centre}

Unfortunately, the author was not able to prove the corresponding
result to Conjecture~\ref{conj:3-1-equiv-Euler} in the case where
$\Ge$ is irreducible but contains $Z_3$.  A method is suggested in
Section~\ref{sec:3-1:centre:2}, but requires further work.

If it the conjecture can be proved in all cases, the work above
implies that $\Term(\C^{4}/G)$ is true for all groups of type $(3,1)$
and $(2,1,1)$.
The results of the next section imply that in that case,
$\Term(\C^{4}/G)$ is true for all irreducible $G$ which contain
$Z_{4}$.

\section{Blowing up in the presence of the centre $Z_n$}
\label{sec:blowing-up}

\subsection{Invariant sets in projective space}

First, a lemma about the Euler number of a invariant sets in
projective space.

\begin{lemma}
  \label{lemma:euler-proj}
  Let $H$ be a finite abelian group acting linearly on $\PP^n$.  Then 
$$\chi((\PP^n)^H) = \chi(\PP^n)=n+1.$$
\end{lemma}
\begin{proof}
  Suppose $H$ has order $r$.  Diagonalise the action of $H$, and order
  the weights of the action so that they form a non-decreasing
  sequence of elements of $\{0,\dots,r-1\}$.  The sequence will
  consist of $d_1$ occurrences of the smallest weight $w_1$, followed
  by $d_2$ occurrences of the second smallest weight $w_2$, and so on,
  ending with $d_s$ occurrences of the greatest weight $w_s$.  Since
  there are $n+1$ (not necessarily distinct) weights in the sequence,
  the multiplicities $d_i$ sum to $n+1$.

  Computing the invariant part of $\PP^n$ with respect to the
  $H$-action, one sees that it consists of a disjoint union over all
  $i\in\{1,\dots,s\}$ of projective spaces $\PP^{d_i-1}$.  Taking the
  sum of the Euler numbers of the invariant components, and
  using the fact that the Euler number of $\PP^d$ is $d+1$, one
  obtains the value $\sum d_i$, which by the previous paragraph indeed 
  coincides with the Euler number $n+1$ of $\PP^n$.
\end{proof}

\subsection{Blowing up the origin}
\label{sec:blowing-up:blowing-up}

Now let $V= \C^n$ and let $\text{Bl}_0V$ be the blow-up of $V$ at the
origin. This has a natural $G$-action and one  has the following
commutative diagram
$$\begin{CD}
         &     \text{Bl}_0V     @>\sigma_0>>      & V  \\
         & \downarrow        &    &         & \downarrow \\
         &    \text{Bl}_0V/G   @>{\sigma'_0}>> & V/G .
\end{CD}$$

A standard discrepancy calculation yields the following result.
\begin{lemma}
  The morphism $\text{Bl}_0V/G \to V/G$ is crepant if and
  only if $G$ contains $Z_n=\<<\omega_n>>$.
\end{lemma}

The following lemma contains the basic idea to constructing
Euler blow-ups.

\begin{lemma}
  \label{lemma:bV}
  Assume that $G$ contains $Z_n$ and write $\bG:=G/Z_n$.  Then $\bV := 
  \text{Bl}_0V/Z_n$ is smooth, and $\bV/\bG \to V/G$ is a projective 
  Euler blow-up.
\end{lemma}
\begin{proof}
  The only place where singularities of $\bV$ could arise is on the
  image $\bE=E/Z_n$ of the exceptional divisor $E$ of $\sigma_0$.
  Identifying $E$ with $\PP(V)$, a local chart for $\text{Bl}_0V$ at
  a point $\xi\in E$ is given in suitable local coordinates by
$$(x_1,\frac{x_2}{x_1},\dots,\frac{x_n}{x_1}),$$
 so therefore $\omega_n$ acts there as $(\omega_n, 1,1,1)$, \ie  as a 
 pseudo-reflection.

The Euler number computation goes as follows:
\begin{align*}
\chiorb(\bV/\bG)
  &=\sum_{[\bg]\in\Cl(\bG)}\chi(\PP(V)^{\bg}/\Cent^{\bG}_{\bg}), \\
  & \qquad\text{since $\bV\sim E\cong\PP(V)$\ \ (homotopy)}\\
  &=\sum_{[\bg]\in\Cl(\bG)} \chi(\PP(V)), \\
  & \qquad\text{averaging and using Lemma~\ref{lemma:euler-proj}},\\
  &=|\Cl(\bG)| n = |\Cl(G)|\\
  &=\chiorb(V/G).
\end{align*}
\end{proof}

Thus if $G$ contains $Z_n$ then $\Term(V/G)\equiv\Term(\bV/\bG)$.

\subsection{Patching}
\label{sec:blowing-up:patching}

In many cases, an Euler blow-up is constructed by patching together
local Euler blow-ups.  The following lemma summarises the necessary
conditions to carry this out this procedure.

\begin{lemma}[Patching Lemma]
  \label{lemma:patching}
  Let $G$ be a finite subgroup of $\SL(n)$ which contains $Z_n$, so
  that $Y:=\bV/\bG\to V/G$ is an Euler blow-up.  Let $\bE:\equiv
  E/Z_n$, where $E$ is the exceptional divisor of $\text{Bl}_0 V\to V$
  and denote by $p\colon\bV\to\bE$ the projection.
  
  There exists a finite collection of points $y\in\bE/\bG$ and 
  corresponding analytic neighbourhoods $\bE_\xi\subset\bE/\bG$ such that
  $Y$ is covered by $\{Y_y\}$, where $Y_y:=p^{-1}(\bE_y)$.
  
  Suppose that for each $y$, there exists an Euler blow-up
  $\varphi_y\colon X_y\to Y_y$, such that if $y\neq y'$, one
  nevertheless has
\begin{equation}
  \label{eq:patching}
 X_{y|Y_y\cap Y_{y'}} =  X_{y'|Y_y\cap Y_{y'}}.
\end{equation}
Then the analytic \cGo\ $X$ obtained by gluing together all the
$ X_y$ is an Euler blow-up of $Y$.
\end{lemma}
\begin{proof}
  The existence of the finite covering $\{\bE_y\}$ follows because
  $\bE$ is compact.  Since, all the $X_y$ are birational to each other 
  above the overlaps,
  equation~\eqref{eq:patching}  implies  that $X$ is well-defined.  
  Furthermore, since no crepant
  divisors are introduced during the local blow-ups and since the  orbifold
  Euler numbers of $X_y$ and $Y_y$ are the same, $X$ is crepant
  over $Y$ and has the same orbifold Euler number.
\end{proof}

\section{Irreducible $G$ which contain $Z_4$}
\label{sec:4-centre}

\subsection{Notation}
\label{sec:4-centre:notation}

When $G$ contains $Z_4$, $\bV/\bG\to V/G$ is an Euler blow-up. 
Let $\bE$ be the exceptional divisor of $\bV\to V/Z_4$, and let
$p\colon \bV\to \bE$ be the projection.  Let $\xi\in \bE$ be a point
in the base, and consider the tangent space of $\bV$ at $\xi$.  This
decomposes into $\bG$-modules 
$$V_\xi^1\oplus V_\xi^2,$$
where $V^1_\xi$ is the tangent space to
$\bE$ and $V^2_\xi$ is the line tangent to the fibre of $p$, and
stabilized by $\bG$.  Let $\xi'\in V_\xi^2$ be the generic point, so
that its stabiliser $\bG_{\xi'}$ is a subgroup of $\SL(3)$.

\subsection{Local Blow-ups}
\label{sec:4-centre:blowups}

Let $\xi\in \bV$ and let $\bar\xi$ denote its image in $Y=\bV/\bG$.
A local analytic neighbourhood of $\bar\xi$ is isomorphic to
$$Y_{\bar\xi} := (V_\xi^1\oplus V_\xi^2)/\bG_\xi= (V_\xi^1/\bG_{\xi'}\oplus
V_\xi^2)/C_{\xi},$$ where $C_\xi$ denotes the  quotient $C_\xi:=
\bG_{\xi}/\bG_{\xi'}$, which is cyclic, being a naturally a subgroup
of $\GL(V_\xi^2)\cong\CX$.

Restricting attention to the ``base'' $\bE$, one has, corresponding to
each $\bar\xi$ in $\bE $, a quotient singularity $\bE_{\bar\xi} =
V_\xi^1/\bG_{\xi'}$ which is an $SL(3)$-singularity.  Since $\bE/\bG$
is compact, the choice of a finite number of points $\bar\xi$ is
sufficient for $\bigcup_{\bar\xi} Y_{\bar\xi}$ to cover the whole of
$\bV/\bG$.

\subsection{Gluing}
\label{sec:4-centre:gluing}

\begin{lemma} 
\label{lemma:4-centre-gluing} 
Let $W_\xi^1\to V_\xi^1/\bG_{\xi'}$, be the
$\bG_\xi$-equivariant Euler resolution such as that in 
Conjecture~\ref{conj:3-1-equiv-Euler}.  Define $$
X_\xi := (W_\xi^1\times V_\xi^2)/\bG_{\xi}.$$
Then $\varphi_\xi\colon X_\xi\to
Y_\xi$ are Euler blowups which have only cyclic quotient
singularities, and for each $\xi,\xi'$, $ X_\xi$ and $ X_{\xi'}$ agree
above the inverse image of $Y_\xi\cap Y_{\xi'}$.  Thus by
Lemma~\ref{lemma:patching}, they glue to form a complex analytic Euler
blow-up $Y \to \bV/\bG$ which has only cyclic quotient singularities.
\end{lemma}
\begin{proof}
  One must check that $\chiorb(X_\xi)=\chiorb(Y_\xi)$ and that , on
  $\varphi_\xi^{-1}(Y_\xi\cap Y_{\xi'})$ the blow-ups corresponding to
  $\xi$ and $\xi'$ agree.  The orbifold Euler number equality is
  checked in
  Proposition~\ref{prop:3-1-invariant-resolution-gives-terminalisation}.
  For the agreement of the blow-ups on the overlaps, knowledge of the
  $SL(3)$ singularities implies that these overlaps only occur over
  curves of 2-dimensional singularities (the components of the
  non-isolated singularities are all curves for $\SL(3)$).  Over these, the
  resolutions which are being glued-in are trivial families of minimal
  resolutions: they are therefore unique, and so resolutions coming
  from neighbourhoods corresponding to different $\xi$'s will agree.
\end{proof}

This gives a crepant analytic blow-up $Y$ which is locally
analytically isomorphic to a cyclic quotient (and hence locally
analytically $\Q$-factorial).

\subsection{Terminalisation and the Orbifold Euler Number}
\label{sec:4-centre:terminalisation}

Since the orbifold Euler number can be calculated by summing the
contributions of the various analytic neighbourhoods, the equality
$$\chiorb(X)=\chiorb(\bV/\bG)$$
will follow by showing that the
resolutions glued in above preserve the orbifold Euler number.  This
is proved in Proposition~\ref{prop:3-1-nocentre-equiv-Euler}.  Thus
$X$ is an Euler blow-up with only cyclic singularities.

Applying the minimal model program to $X$
(Theorem~\ref{thm:toroidal-mmp}), one obtains a crepant terminalisation
$t\colon Z\to X$ which satisfies $\chiorb(Z)=\chiorb(X)=\chiorb(V/G)$,
and has only toric singularities.

\begin{rmk}
  By studying which toric flips can occur, one might be able to prove 
  that the singularities of $T$ are in fact at most cyclic.  They 
  would then have to necessarily be isolated.  For if a 
  4--dimensional Gorenstein cyclic singularity consisted of a curve 
  of singularities, these would also have to be (3-dimensional) 
  terminal Gorenstein cyclic quotients.  But the classification of 
  3-dimensional terminal cyclic quotients~\cite{mor_stev:terminal} 
  shows that they are all of the form $\frac{1}{r}(1,-1,a)$, and so 
  can only be Gorenstein if they are smooth.
\end{rmk}

\subsection{Case where $G$ is of type $(3,1)$ revisited}
\label{sec:3-1:centre:2}

A method similar to the one above can be applied to the case treated
in Section~\ref{sec:3-1:centre}, namely the case where $G$ is a group
of type (3,1) and the stabiliser group $\Ge<\SL(3)$ is irreducible and
contains $Z_3$.  This goes some way towards a solution of
Conjecture~\ref{conj:3-1-equiv-Euler}.

Since $\Ge >Z_3$, there exists an Euler blow-up $\bV^1/\bGe\to 
V^1/\Ge$ and this is equivariant under the $G$ (and hence $G'$) 
action, since it is obtained from blowing up the origin of $V^1$, 
which is of course fixed by $G$.  Its singularities are of the type 
$(2,1)$ and $(1,1,1)$ and the singular locus is 
invariant under $G$.
  
An analytic resolution of $\bV^1/\bGe$ can be constructed by the same 
gluing procedure as in Lemma~\ref{lemma:4-centre-gluing}: in any 
local analytic neighbourhood of $\bV^1/\bGe$, construct an Euler 
resolution, doing this equivariantly under the $G'$-action.  These 
glue together, since they can only intersect over smooth points.  
This gives a resolution $W^1\to\bV^1/\bGe$.

It remains to show that it admits a $\bG$-action and that it 
satisfies the Euler number property of equation~\eqref{eq:chi-pi}.

\section{Groups $G$ of type (2,2) and (2,1,1)}
\label{sec:2-}

\subsection{Groups of type $(2,2)$}
\label{sec:2-2}

Let $\C^4=V^1\oplus V^2$ with $V^i$ 2-dimensional irreducible
$G$-modules and denote by $\eta_i$ the generic point of $V^i$.  If
$G<\SL(2)\times\SL(2)$ then the following lemma constrains the
stabilisers of $\eta_i$.

\begin{lemma}
  \label{lemma:terminal}
  Suppose that $V=V^1\oplus V^2$ and that
  $G<\SL(V^1)\times\SL(V^2)<\SL(V)$.  Denote by $G_{i}$ the stabiliser
  of the generic point $\eta_i\in V^i$ for $i=1,2$.  If both 
  stabilisers $G_1$ and $G_2$ are trivial, then $V/G$ must be terminal.
\end{lemma}
\begin{proof} 
  This can be proved by an easy discrepancy calculation, or
  equivalently, by using the concept of ``weights'' for the group
  action~\cite{reid_ito} as follows.  The number of crepant divisors
  of $V/G$ is equal to the ``number of elements of $G(1)$ of weight
  one''.  Here, $G(1):=\Hom(\mu_r, G)$, where $r$ is the least common
  multiple of the orders of the elements of $G$.  The {\em weight\/}
  $\wt(\hat g)$ of $\hat g\in G(1)$ is defined by evaluating $\hat g$
  on a primitive generator $\epsilon$ of $\mu_r$, diagonalising the
  resulting matrix $\hat g(\epsilon)$ and expressing the diagonal
  elements in terms of powers of $\epsilon$ ranging between $0$ and
  $r-1$.  Because $G<\SL(V)$, the sum of these powers divided by $r$
  is a non-negative integer, called the {\em weight\/} of $\hat g$.
  Note that if we simply want to calculate the {\em number\/} of
  elements of $G(1)$ of a given weight, we can identify $G$ with
  $G(1)$ by fixing a primitive generator of $\mu_r$, and pretend to be
  calculating the weights of the elements of $G$.
  
  Suppose $V/G$ is not terminal, so that there exists and element
  $g\in G$ of weight one.  For each $i=1,2$, denote by $g_i$ the
  part of the matrix of $g$ which represents its action on the module
  $V^i$.  Since $\wt(g)=\wt(g_1)+\wt(g_2)$, one of the $g_j$'s must be
  equal to the identity (whereas the other one must be a non-trivial
  matrix).  But $G_j=1$ and $g_j=1$ for some $j$ would imply that
  $g=1$, so one of the two stabilisers must be non-trivial.
\end{proof}

Thus, $G<\SL(2)\times\SL(2)$ then $V/G$ not terminal implies
(Lemma~\ref{lemma:terminal}) that one of $G_{\eta_i}$ must be
non-trivial, say $G_\eta=G_{\eta_2}\neq 1$.  Denoting by $W^1\to
V^1/G_\eta$ the minimal resolution, one obtains a crepant blow-up
$(W^1\times V^2)/(G/G_\eta)\to V/G$ with singularities of type
$(2,1,1)$ or $(1,1,1,1)$.

In order to prove that the orbifold Euler number remains unchanged 
under the blowup $(W^1\times V^2)/G\to V/G$, one must prove a formula 
similar to that of equation~\eqref{eq:chi-pi}, namely:
$$\chi((W^{1})^{d}/Z_d)=|\pi_*^{-1}([d])|,$$
for all $d\in D:=G/G_{\eta}<\GL(V^1)$.  

This could presumably be achieved by examining the finite $\SL(2)$ and
$\GL(2)$ subgroups and determining how many fit into an exact sequence 
of the form given in equation~\eqref{eq:3-1-nocentre-exact}, with the
cyclic group $C$ replaced by the $\GL(2)$ subgroup $D$.

\begin{question}
  What happens if $G$ is not a subgroup of $\SL(2)\times\SL(2)$?
\end{question}

\subsection{Groups of type $(2,1,1)$}
\label{sec:2-1-1}

Suppose the irreducible decomposition of $\C^4$ is $V^1\oplus V^2
\oplus V^3$ with $\dim V^1=2$ and $\dim V^i=2$ for $i=2,3$.  Denote by
$\eta_i$ the generic point of $V^i$,  by $G_3$ the stabiliser of
$\eta_3$ in $G$ and by $G_{32}$ the stabiliser of $\eta_2$ in $G_3$.
Then we have exact sequences
\begin{alignat}{3}
  \label{eq:2-1-1-exact}
    1 &\to G_3 &\to G &\xrightarrow{\pi_3} C \to 1,\\
\intertext{and}
   \label{eq:2-1-1-exact2}
    1 &\to G_{32} &\to G_3 &\xrightarrow{\pi_{32}} D \to 1,
\end{alignat}
with $C$ and $D$ cyclic subgroups, $G_{32}<\SL(2)$ and $G_3<\SL(3)$.

The analytic germ of $V/G$ is isomorphic to
$$(V^1\times V^2 \times V^3)/G \cong \bigl( (V^1/G_{32}\times V^2)/D\times
V^3\bigr )/C.$$ 

The term $V^1/G_{32}$ has a minimal resolution $Z^1$ which is unique
and therefore admits an action of $G_3$.  The germ $(Z^1\times V^2)/D$
is a cyclic $\SL(3)$ singularity and has an Euler resolution $W^1$.
It remains to be shown that $W^1$ admits a $G$-action and that
$\chiorb((W^1\times V^3)/C)=\chiorb(V/G)$.  As in
Section~\ref{sec:3-1}, the later statement would follow
from the equalities
$$\chi((W^1)^{c'}/C')=|\pi_{3*}^{\prime -1}([c'])|,$$
for $c'\in C'$, where the primed objects are defined similarly to those 
in Section~\ref{sec:4-centre:notation}.

\providecommand{\bysame}{\leavevmode\hbox to3em{\hrulefill}\thinspace}

\end{document}